\documentstyle{article}

\newcommand{\R}{{\bf R}}
\newcommand{\id}{{\rm Id}}

\newtheorem{thm}{Theorem}[section]

\newtheorem{lem}{Lemma}[section]
\newtheorem{cor}{Corollary}[section]
\newtheorem{definition}{Definition}[section]

\def\thebibliography#1{\section*{Bibliography}\list
 {[\arabic{enumi}]}{\settowidth\labelwidth{[#1]}\leftmargin\labelwidth
 \advance\leftmargin\labelsep
 \usecounter{enumi}}
 \def\newblock{\hskip .11em plus .33em minus .07em}
 \sloppy\clubpenalty4000\widowpenalty4000
 \sfcode`\.=1000\relax}

\def\resume{\if@twocolumn
\section*{Abstract}
\else \small 
\begin{center}
{\bf R\'esum\'e\vspace{-.5em}\vspace{0pt}} 
\end{center}
\quotation 
\fi}
\def\endresume{\if@twocolumn\else\endquotation\fi}

\newenvironment{preuve}{\par\noindent{\bf Proof~:}}{\hfill\stopthm} 
\newcommand{\stopthm}{\hfill$\diamond$\par\smallskip} 
\author{St\'ephane Grognet}
\title{Weakly convex closed subsets of spaces with bounded nonpositive curvature}
\begin{document}

\maketitle

\vskip 30 pt
\noindent Universit\'e de Nantes, Laboratoire Jean Leray, U. M. R. 6629, 2, rue
de la Houssini\`ere, BP 92208, F-44322 Nantes cedex 03.
 
\noindent Stephane.Grognet@univ-nantes.fr

\vskip 10 pt

R\'ESUM\'E. - Sur~$\R^n$ muni d'une m\'etrique riemannienne \`a courbure n\'egative born\'ee,
les ferm\'es faiblement convexes sont topologiquement triviaux.
Leur stabilit\'e par intersection caract\'erise les espaces euclidiens.
\vskip 10 pt

ABSTRACT. - On~$\R^n$ endowed with a riemannian metric of bounded nonpositive curvature,
the weakly convex closed subsets are topologically trivial.
The stability of such subsets under intersection characterizes the euclidean spaces.

\vskip 8 pt

AMS-MSC~: 53C20, 52A55

{\em keywords}~: weak convexity, geodesic flow, negative curvature, horospheres

%%%%%%%%%%%%%%%%%%%%%%%%%%%%%%%%%%%%%%%%%%%%%%%%%%%%%%%%%%%%%%%%%%%%%%%%%%%%%%%
\section{Introduction}\label{sec:intro}
%%%%%%%%%%%%%%%%%%%%%%%%%%%%%%%%%%%%%%%%%%%%%%%%%%%%%%%%%%%%%%%%%%%%%%%%%%%%%%%

The convexity of closed subsets in the euclidean spaces and in the riemannian manifolds has been subject of many studies
among which we will only cite~\cite{CH,W1,W2}, the third presenting an important bibliography.
In presence of negative curvature, there appears a notion of weak convexity different from the geodesic convexity~\cite{Gro3,Gro4}.
Horospheres are instances of closed subsets which are not geodesically convex but are weakly convex~;
others examples in dimension two can be found in~\cite{Gro4}.

In the whole paper we consider a complete, connected, simply connected manifold~$M$ of~${\cal C}^\infty $-class,
endowed with a riemannian metric with bounded nonpositive sectional curvature~:~$-k^2\leq K\leq 0$.
According to the Cartan-Hadamard theorem, the manifold~$M$ is diffeomorphic to~$\R^n$.

We present here two expected properties of the weakly convex closed subsets~: their topology is trivial,
and their bad behaviour with intersection is characteristic.

Given a closed subset of~$M$, it is said that {\em the projection is unique} 
when the distance to~$G$ of each point of the complementary~$G^C$ is realized by a unique point of~$G$.
The projection from~$G^C$ to~$\partial G$, which to a point~$x$ 
associates the unique point~$y$ such that~$d(x,G)=d(x,y)$ is written down~$\pi $.
\begin{definition} \cite{Gro3}
A closed subset~$G$ of~$M$ is said {\em weakly convex} if the projection 
is unique, and if for any point~$x$ of the complementary~$G^C$, the open horobowl 
associated to the geodesic ray~$[\pi (x),x)$ does not intersect the closed subset~$G$~: 
the half geodesic~$[\pi (x),x)$ projects itself on the point~$\pi (x)$.
\end{definition}
Given a closed subset of~$M$, the Motzkin theorem in bounded nonpositive sectional curvature~\cite{Gro3}
ensures the equivalence between the three properties of weak convexity, uniqueness of projection,
and differentiability of the distance to~$G$ on the complementary~$G^C$.

The condition of geodesic convexity is stronger than this weak convexity~\cite{Gro3}~;
weak and strong convexity coincide in euclidean spaces.

We prove that nonempty weakly convex closed subsets are deformation retracts of~$M$ (theorem~\ref{thm:retracte}).
In particular they are connected, simply connected as expected (corollary~\ref{cor:simpleConnexe}).
We deduce from it that the property of stability of the weakly convex closed subsets under finite intersection
characterizes the euclidean spaces (theorem~\ref{thm:inter}).

%%%%%%%%%%%%%%%%%%%%%%%%%%%%%%%%%%%%%%%%%%%%%%%%%%%%%%%%%%%%%%%%%%%%%%%%%%%%%%%
\section{Topological viewpoint}
%%%%%%%%%%%%%%%%%%%%%%%%%%%%%%%%%%%%%%%%%%%%%%%%%%%%%%%%%%%%%%%%%%%%%%%%%%%%%%%

\begin{lem}\label{thm:homotopie}
Let~$M $ be a~${\cal C}^\infty $-manifold diffeomorphic to~$\R^n$ equipped with a riemannian metric of bounded nonpositive sectional curvature.
For~$x,y\in M$ and~$t\in[0,1]$, let us define~$(1-t)x+ty$ as a barycentre on the geodesic segment~$[x,y]$.
Let~$G$ be a nonempty weakly convex closed subset of~$M$. The mappings
\[
\matrix{
H: & M \times [0,1] & \to     & G \cr
{} & (x,t)          & \mapsto & x \hbox{ \rm if } x\in G \cr 
{} & (x,t)          & \mapsto & (1-t)x+t\pi (x) \hbox{ \rm if } x\in G^C \cr
}
\]
and
\[
\matrix{
P: & M        & \to     & G \cr
{} & x\in G   & \mapsto & x \cr
{} & x\in G^C & \mapsto & \pi (x) \cr
}
\]
are continuous.
\end{lem}
\begin{preuve}
The continuity of~$P$ implies that of~$H$.
It is clear that~$P$ is continuous at every point of the complementary~$G^C$ and the interior~$\displaystyle { \mathop{G}^\circ }$.
Let~$x$ be a point of~$\partial G$. We have~$P(x)=x$. Let~$\varepsilon >0$.
Let~$y$ be a point of~$B(x,\varepsilon ) \cap G$. We have~$P(y)=y \in B(x,\varepsilon ) $.
Let~$y$ be a point of~$B(x,\varepsilon ) \cap G^C$. We have~$P(y)=\pi (y)$ with
\[
d \left( y,\pi (y) \right) \leq d(y,x) < \varepsilon
\]
thus
\[
d \left( \pi (y) ,x \right) = d \left( P(y),P(x) \right) \leq d \left( \pi (y),y \right) + d(y,x) < 2 \varepsilon .
\]
Therefore~$P$ is continuous at~$x$.
\end{preuve}
\begin{thm}\label{thm:retracte}
Let~$M $ be a~${\cal C}^\infty $-manifold diffeomorphic to~$\R^n$ equipped with a riemannian metric of bounded nonpositive sectional curvature.
Every nonempty weakly convex closed subset~$G$ of~$M$ is a deformation retract of~$M$.
\end{thm}
\begin{preuve}
This results from the previous lemma, since~$H(.,0)=\id _M $.
\end{preuve}

As an immediate consequence, we get~:
\begin{cor}\label{cor:simpleConnexe}
Let~$M $ be a~${\cal C}^\infty $-manifold diffeomorphic to~$\R^n$ equipped with a riemannian metric of bounded nonpositive sectional curvature.
Every nonempty weakly convex closed subset~$G$ of~$M$ is connected, simply connected.
\end{cor}

%%%%%%%%%%%%%%%%%%%%%%%%%%%%%%%%%%%%%%%%%%%%%%%%%%%%%%%%%%%%%%%%%%%%%%%%%%%%%%%
\section{Geometrical viewpoint}
%%%%%%%%%%%%%%%%%%%%%%%%%%%%%%%%%%%%%%%%%%%%%%%%%%%%%%%%%%%%%%%%%%%%%%%%%%%%%%%

The behaviour with respect to intersection of the weakly convex closed subsets is different from that of the geodesically convex ones.
\begin {definition}
The property of {\em stability under finite intersection} for the weakly convex closed subsets is the following~:
for all weakly convex closed subsets~$G_1,G_2$, the intersection~$G_1\cap G_2$ is weakly convex.
\end{definition}
In general, this property is not true.
\begin{thm}\label{thm:inter}
Let~$M $ be a~${\cal C}^\infty $-manifold diffeomorphic to~$\R^n$ equipped with a riemannian metric of bounded nonpositive sectional curvature.
The property of stability under finite intersection for the weakly convex closed subsets is true if and only if
the metric is euclidean.
\end{thm}
\begin{preuve}
For euclidean spaces, the condition of weak convexity coincides with the geodesic convexity, so the property is true.
In order to prove the other implication, we suppose the property to be true.
If there exists a point~$x\in M$ and a~$2$-plane~$\Pi\subset T_xM$ such that~$K(\Pi )<0$, 
let~$v$ be a unit tangent vector of~$T^1_xM$, belonging to~$\Pi $.
Let~$({\cal H}_i)_{i=1,\dots ,n-2}$ be a collection of horospheres tangent to~$\Pi $,
such that their intersection admits~$\Pi $ as tangent space at~$x$.
The horospheres are weakly convex.
Let~$G_1$ be the complementary of the open stable horobowl associated with~$v$.
The geodesic flow on~$T^1M$ is written~$\varphi_t$.
For some~$\varepsilon >0$, let~$G_2$ be the complementary of the open unstable horobowl associated with~$\varphi_\varepsilon v$.
The closed subsets~$G_1$ and~$G_2$ are weakly convex~\cite{Gro3}.
Because of the property, the intersection
\[
I=G_1 \cap G_2 \cap \left( \mathop{\cap }_{i=1}^{n-2} {\cal H}_i \right) 
\]
is weakly convex. But because of the condition~$K(\Pi )<0$, for~$\varepsilon $ small enough, 
the second fondamental form of the horospheres~$\partial G_1$,~$\partial G_2$ 
are nonzero in the orthogonal of~$v$ in~$\Pi $.
Thus, for~$\varepsilon $ small enough, the intersection~$I$ is nonempty and not connected,
which contradicts the corollary~\ref{cor:simpleConnexe}.
Therefore the curvature is always zero.
\end{preuve}

%%%%%%%%%%%%%%%%%%%%%%%%%%%%%%%%%%%%%%%%%%%%%%%%%%%%%%%%%%%%%%%%%%%%%%%%%%%%%%%
%%%%%%%%%%%%%%%%%%%%%%%%%%%%%%%%%%%%%%%%%%%%%%%%%%%%%%%%%%%%%%%%%%%%%%%%%%%%%%%

\end{document}